\documentclass{ifacconf}

\usepackage{amsmath}
\usepackage{amssymb}

\usepackage{algorithm}
\usepackage{algpseudocode}

\usepackage{graphicx}     
\graphicspath{ {figure/} }
\usepackage{tabularx}

\usepackage{natbib}        
\usepackage{url}

\begin{document}

\begin{frontmatter}
\title{A Successive Gap Constraint Linearization Method for Optimal Control Problems with Equilibrium Constraints\thanksref{sponsor and financial support acknowledgment}} 

\thanks[sponsor and financial support acknowledgment]{This work was partly supported by JSPS KAKENHI Grant Number JP22H01510 and JP23K22780. 
The CSC scholarship supported the author Kangyu Lin (No. 201906150138).}

\author[author]{Kangyu Lin} 
\author[author]{Toshiyuki Ohtsuka} 

\address[author]{Department of Systems Science, Graduate School of Informatics, Kyoto University, Kyoto, Japan (e-mail: \{k-lin, ohtsuka\}@sys.i.kyoto-u.ac.jp).}

\begin{abstract}  
In this study, we propose a novel gap-constraint-based reformulation for optimal control problems with equilibrium constraints (OCPECs).
We show that the proposed reformulation generates a new constraint system equivalent to the original one but more concise and with favorable differentiability.
Moreover, constraint regularity can be recovered by a relaxation strategy.
We show that the gap constraint and its gradient can be evaluated efficiently.
We then propose a successive gap constraint linearization method to solve the discretized OCPEC.
We also provide an intuitive geometric interpretation of the gap constraint.
Numerical experiments validate the effectiveness of the proposed reformulation and solution method.
\end{abstract}

\begin{keyword}
Optimization and Model Predictive Control, Optimal Control, Hybrid Systems, Variational Inequalities, Gap Function
\end{keyword}

\end{frontmatter}


\section{Introduction}

\subsection{Background}
Optimal control problems (OCPs) for a class of nonsmooth dynamical systems governed by \textit{differential variational inequalities} (DVIs) are referred to as \textit{optimal control problems with equilibrium constraints} (OCPECs).
Benefiting from the mathematical modeling capability of \textit{variational inequalities} (VIs), DVIs encompass many common nonsmooth systems, such as ordinary differential equations (ODEs) with discontinuous right-hand-side, dynamical complementary systems, and dynamical systems with optimization constraints.
Thus, OCPECs provide a broad paradigm for OCPs in many complex engineering systems, such as the trajectory optimization of a mechanical system with frictional contacts \citep{Posa2014}, bilevel optimal control \citep{dempe2023bilevel}, and game-theoretic planning for autonomous driving \citep{wang2021game}.

The direct method \citep{betts2010practical}, that is, first discretize then optimize, is practical for numerically solving OCPs.
However, its application to OCPECs is extremely challenging. 
On the one hand, the nonsmooth state trajectory limits the system discretization to achieve at most first-order accuracy. 
This leads to a large-scale optimization problem where mixed-integer program solvers often fail due to the large search tree.
On the other hand, the discretized OCPEC violates almost all constraint qualifications; therefore, typically, it cannot be solved directly by large-scale nonlinear programming (NLP) solvers. 
The most common approach to solve the discretized OCPEC is to reformulate it as the \textit{mathematical programming with complementary constraints} (MPCC) and then utilize MPCC-tailored solution methods, such as relaxation methods \citep{hoheisel2013theoretical}.
However, it is well known that MPCC also lacks constraint regularity; thus, this approach still encounters numerous difficulties.
Owing to the fundamental absence of constraint regularity in OCPECs, we cannot expect to construct an OCPEC reformulation satisfying constraint qualifications.
Nonetheless, we can still attempt to construct a new OCPEC reformulation that exhibits improved numerical performance, which is our motivation.

In this study, based on the regularized gap function for VIs, we propose a novel gap-constraint-based reformulation for the OCPEC, which generates a concise differentiable constraint system and can achieve constraint regularity after a suitable relaxation.
To the best of our knowledge, no MPCC or OCPEC reformulations are currently using gap functions.
We believe the main reason is that evaluating gap functions and their gradients is typically quite expensive, which also limits the application of gap-function-based methods for solving VIs (Chapter 10, \cite{facchinei2003finite}).
However, we show that this process can be accelerated by exploiting either VI or OCP structures.
Subsequently, we propose a successive gap constraint linearization method to solve the discretized OCPEC.
We show that the proposed reformulation captures the geometric characteristics of equilibrium constraints in a favorable shape.
Numerical experiments confirm that the proposed method outperforms the MPCC-tailored method.

\subsection{Outline}
This paper is organized as follows:
Section \ref{section: problem formulation} provides the OCPEC formulation.
Section \ref{section: proposed reformulation for OCPEC} and \ref{section: proposed solution method for OCPEC} respectively present the proposed reformulation and solution method.
Section \ref{section: geometric interpretation} provides a geometric interpretation of the proposed reformulation.
Section \ref{section: numerical experiment} illustrates the numerical experiments.
Section \ref{section: conclusion and future work} concludes this study and discusses future works.

\subsection{Notation}
We denote the nonnegative orthant of an Euclidean vector space $\mathbb{R}^n$ by $\mathbb{R}^n_{+}$.
We denote the transposed Jacobian of a vector-valued differentiable function $f(x): \mathbb{R}^n \rightarrow \mathbb{R}^m$ by $\nabla_x f \in \mathbb{R}^{m \times n}$.
We denote the Hessian of a real-valued differentiable function $f(x): \mathbb{R}^n \rightarrow \mathbb{R}$ by $\nabla_{xx} f \in \mathbb{R}^{n \times n}$.
We denote the $A$-norm in $\mathbb{R}^n$ by $\| x \|_A = \sqrt{x^T A x}, \forall x \in \mathbb{R}^n$, with $A \in \mathbb{R}^{n \times n}$ a symmetric positive definite matrix.  

\section{Problem formulation}\label{section: problem formulation}

\subsection{Variational inequalities}
The finite-dimensional VI serves as a unified mathematical formalism for various equilibrium problems  \citep{facchinei2003finite}.
Its definition is provided as follows:
\begin{defn}[Variational Inequalities]\label{definition: VI}
    Let $K \subseteq \mathbb{R}^{n_\lambda}$ be a nonempty closed convex set and $F: \mathbb{R}^{n_\lambda} \rightarrow \mathbb{R}^{n_\lambda}$ be a continuous function; the variational inequalities, denoted by VI$(K, F)$, is to find a vector $\lambda \in K$ such that
    \begin{equation}\label{equation: VI definition}
        (\omega - \lambda)^T F(\lambda) \geq 0, \quad \forall \omega \in K.
    \end{equation}    
\end{defn}
We denote the solution set to the VI$(K, F)$ as SOL$(K, F)$.
The VI includes many common mathematical formalisms:
\begin{itemize}
    \item If $K := \mathbb{R}^{n_\lambda}$, then the VI is reduced to the \textit{system of equations}: $F(\lambda) = 0$;
    \item If $K := \mathbb{R}^{n_\lambda}_{+}$, then the VI is reduced to the \textit{nonlinear complementary problem} (NCP): $0 \leq \lambda \perp F(\lambda) \geq 0$, where $\lambda \perp F(\lambda)$ represents the element-wise complementarity conditions, that is, $\lambda \odot F(\lambda) = 0$ and $\odot$ is the Hadamard product.
\end{itemize}
In this study, we focus on the general case where set $K$ is finitely representable by convex inequalities:
\begin{equation}\label{equation: finitely representable VI set K}
    K:= \{\lambda \in \mathbb{R}^{n_\lambda} | g(\lambda) \geq 0\},
\end{equation}
with $g: \mathbb{R}^{n_\lambda} \rightarrow \mathbb{R}^{n_g}$ a smooth convex function.
The most common approach to solving the VI is to use well-developed solution methods for the NCP.
This is based on that under certain assumptions (Section 1.3, \cite{facchinei2003finite}), SOL$(K, F)$ is equivalent to the solution set to the following Karush--Kuhn--Tucker (KKT) system: 
\begin{subequations}\label{equation: KKT condition based reformulation for VI}
    \begin{align}
        & F(\lambda) - \nabla_{\lambda}g^T(\lambda)\zeta = 0, \\
        & 0 \leq \zeta \perp g(\lambda) \geq 0,
    \end{align}
\end{subequations}
with a constraint multiplier $\zeta \in \mathbb{R}^{n_g}$, and we call (\ref{equation: KKT condition based reformulation for VI}) the \textit{KKT-condition-based reformulation} for the VI.

An alternative approach is to reformulate the VI as a differentiable optimization problem that optimizes a suitable merit function tailored for the VI. 
The definition of the merit function is as follows:
\begin{defn}[Merit function for the VI]
    A merit function for the VI$(K, F)$ on a (closed) set $X \supseteq K$ is a nonnegative function $\theta: X \rightarrow \mathbb{R}_{+}$ such that $\lambda \in$ SOL$(K, F)$ if and only if $\lambda \in X$ and $\theta(\lambda) = 0$, that is, if and only if the solution set SOL$(K, F)$ coincides with the set of global solutions to the optimization problem:
    \begin{equation} \label{equation: merit function based reformulation for VI}
        \min_{y \in X} \ \theta(y),
    \end{equation}
    and the optimal objective value of this problem is zero.  
\end{defn}
We call (\ref{equation: merit function based reformulation for VI}) the \textit{merit-function-based reformulation} for the VI.
Usually, set $X$ is specified as either the VI set $K$ or the entire space $\mathbb{R}^{n_{\lambda}}$, respectively leading to a constrained and an unconstrained optimization problem for solving the VI.

\subsection{Differential variational inequalities}
The DVI was first proposed in \cite{pang2008differential} to model a class of nonsmooth systems under a unified framework. 
It consists of an ODE and a dynamical VI:
\begin{subequations}\label{equation: DVI definition}
    \begin{align}
        & \Dot{x} = f(x(t), u(t), \lambda(t)), & \label{equation: DVI definition ODE} \\
        & \lambda(t) \in \text{SOL}(K, F(x(t), u(t), \lambda(t))),  \label{equation: DVI definition VI}
    \end{align}
\end{subequations}
where $x(t): [0, T] \rightarrow \mathbb{R}^{n_x}$ is the differential state, 
$u(t): [0, T] \rightarrow \mathbb{R}^{n_u}$ is the control input, 
$\lambda(t): [0, T] \rightarrow \mathbb{R}^{n_\lambda}$ is the algebraic variable, 
$f: \mathbb{R}^{n_x} \times \mathbb{R}^{n_u} \times \mathbb{R}^{n_\lambda} \rightarrow \mathbb{R}^{n_x}$ is the ODE right-hand function,
and $\text{SOL}(K, F)$ is the solution set of a (dynamical) VI defined by a set $K \subseteq \mathbb{R}^{n_\lambda}$ and a function $F: \mathbb{R}^{n_x} \times \mathbb{R}^{n_u} \times \mathbb{R}^{n_\lambda} \rightarrow \mathbb{R}^{n_\lambda}$. 
We make the following assumptions on the DVI throughout this study: 
\begin{assum}
    $K$ is a nonempty closed convex set given by (\ref{equation: finitely representable VI set K}); 
    $f, F$ are twice continuously differentiable functions.
\end{assum}
Note that, $\lambda(t)$ may be nonsmooth, discontinuous, or even unbounded w.r.t. time due to the time-varying and multi-valued set $\text{SOL}(K, F)$ in (\ref{equation: DVI definition VI}). 
Thus, substituting $\lambda(t)$ into ODE (\ref{equation: DVI definition ODE}) may induce nonsmoothness to $x(t)$.
In this study, we discretize DVI (\ref{equation: DVI definition}) by the time-stepping method proposed in \cite{pang2008differential}. Specifically, ODE (\ref{equation: DVI definition ODE}) is discretized by the implicit Euler method, and VI (\ref{equation: DVI definition VI}) is enforced at each time point $t_n \in [0, T]$:
\begin{subequations}\label{equation: discretized DVI}
    \begin{align}
        & x_n = x_{n-1} + f(x_n, u_n, \lambda_n)\Delta t,  \label{equation: discretized DVI ODE} \\
        & \lambda_n \in \text{SOL}(K, F(x_n, u_n, \lambda_n)),  \label{equation: discretized DVI VI} 
    \end{align}
\end{subequations}
for $n = 1,2,..., N$, where $x_n\in \mathbb{R}^{n_x}$ and $\lambda_n\in \mathbb{R}^{n_\lambda}$ are respectively the value of $x(t)$ and $\lambda(t)$ at the time point $t_n$, $u_n \in \mathbb{R}^{n_{u}}$ is the piecewise constant approximation of $u(t)$ in the interval $(t_{n-1}, t_n]$, $N$ is the number of stages, and $\Delta t : = T/N$ is the time step.

\subsection{Optimal control problem with equilibrium constraints}
We consider the continuous-time OCPEC in the form of
\begin{subequations}\label{equation: continuous OCPEC}
    \begin{align}
        \min_{x(\cdot), u(\cdot),  \lambda(\cdot)} \ &  L_{T}(x(T)) +  \int_0^T L_{S}(x(t), u(t), \lambda(t)) dt \\
        \text{s.t.} \  &  G(x(t), u(t)) \geq 0,  \\
                       &  C(x(t), u(t)) = 0,  \\
                       &  \Dot{x} = f(x(t), u(t), \lambda(t)), \label{equation: continuous OCPEC ODE} \\
                       &  \lambda(t) \in \text{SOL}(K, F(x(t), u(t), \lambda(t))) \label{equation: continuous OCPEC VI}
    \end{align}
\end{subequations}
for $t \in [0, T]$, 
where: $L_T: \mathbb{R}^{n_x} \rightarrow \mathbb{R}$ and $L_S: \mathbb{R}^{n_x} \times \mathbb{R}^{n_u} \times \mathbb{R}^{n_\lambda} \rightarrow \mathbb{R}$ are respectively the terminal and stage cost function, 
$G: \mathbb{R}^{n_x} \times \mathbb{R}^{n_u}  \rightarrow \mathbb{R}^{n_G}$ and $C: \mathbb{R}^{n_x} \times \mathbb{R}^{n_u}  \rightarrow \mathbb{R}^{n_C}$ are respectively the inequality and equality path constraint,
and (\ref{equation: continuous OCPEC ODE}), (\ref{equation: continuous OCPEC VI}) represent the DVI defined in (\ref{equation: DVI definition}).
We make the following differentiability assumption: 
\begin{assum}
    Functions $L_T, L_S, G, C$ are twice continuously differentiable.
\end{assum}

Solving OCPEC (\ref{equation: continuous OCPEC}) is extremely difficult, as equilibrium constraints (\ref{equation: continuous OCPEC VI}) lack constraint regularity. 
Therefore, dealing with equilibrium constraints (\ref{equation: continuous OCPEC VI}) is a focal point of this study.
By using the direct multiple shooting method \citep{betts2010practical} with the DVI discretization (\ref{equation: discretized DVI}), we discretize the continuous-time OCPEC (\ref{equation: continuous OCPEC}) into an NLP problem with a given initial state $x_0$ as follows:
\begin{subequations}\label{equation: discretized OCPEC}
    \begin{align}
    \min_{\{ x_n, u_n,\lambda_n\}^N_{n=1}} & L_T(x_N)  + \sum^{N}_{n=1} \underbrace{L_{S}(x_n, u_n, \lambda_n) \Delta t}_{:= L(x_n, u_n, \lambda_n)} , \\
     \text{s.t.}     \quad              & G(x_n, u_n) \geq 0, \\
                                        & C(x_n, u_n) = 0, \\
                                        & x_{n-1} + \underbrace{f(x_n, u_n, \lambda_n) \Delta t - x_n}_{:= \mathcal{F}(x_n, u_n, \lambda_n)}  = 0, \label{equation: discretized OCPEC ODE}\\
                                        & \lambda_n \in \text{SOL}(K, F(x_n, u_n, \lambda_n)), \quad n = 1, \dots ,N. \label{equation: discretized OCPEC VI}
    \end{align}
\end{subequations}
The limits of time-stepping methods for OCPs in nonsmooth systems have been discussed in \cite{nurkanovic2020limits}, which mainly lie in the low discretization accuracy and wrong sensitivity.
To mitigate these numerical issues, one practical approach is to specify a sufficiently small $\Delta t$, smoothing the nonsmooth system, and then solve the smoothed OCP using a continuation method that tracks the smoothing parameter.
In the case of the OCPEC, a common implementation of this approach is to replace VI (\ref{equation: discretized OCPEC VI}) with its KKT-condition-based reformulation (\ref{equation: KKT condition based reformulation for VI}) and further reformulate the complementary constraints using a MPCC-tailored relaxation strategy.
For example, employing Scholtes' relaxation strategy \citep{scholtes2001convergence} generates the following \textit{relaxed KKT-condition-based} problem:
\begin{subequations}\label{equation: discretized OCPEC relax KKT}
    \begin{align}
    \min_{\{ x_n, u_n,\lambda_n\}^N_{n=1}} & L_T(x_N) + \sum^{N}_{n=1}L(x_n, u_n, \lambda_n),  \\
     \text{s.t.}     \quad              & G(x_n, u_n) \geq 0 ,\\
                                        & C(x_n, u_n) = 0 ,\\
                                        & x_{n-1} + \mathcal{F}(x_n, u_n, \lambda_n) = 0, \\
                                        & F(x_n, u_n, \lambda_n) - \nabla_{\lambda}g^T(\lambda_n)\zeta_n = 0, \label{equation: discretized OCPEC relax KKT VI stationary}\\
                                        & \zeta_n \geq 0, \label{equation: discretized OCPEC relax KKT VI dual} \\ 
                                        & g(\lambda_n) \geq 0, \label{equation: discretized OCPEC relax KKT VI set} \\ 
                                        & s I_{n_g \times 1} - \zeta_{n} \odot g(\lambda_n) \geq 0,  \label{equation: discretized OCPEC relax KKT VI complementary}
    \end{align}
\end{subequations}
with $n = 1, \dots ,N$, $s \geq 0$ the relaxation parameter, and $\zeta_n \in \mathbb{R}^{n_g}$ the constraint multiplier for $g(\lambda_n) \geq 0$.
Therefore, solutions to the discretized OCPEC (\ref{equation: discretized OCPEC}) can be obtained by the continuation method, that is, solving a sequence of relaxed problems (\ref{equation: discretized OCPEC relax KKT}) with $s \rightarrow 0$.
However, numerical difficulties arise when $s$ is close to zero, as the relaxed problem (\ref{equation: discretized OCPEC relax KKT}) degenerates to an MPCC when $s = 0$.
This motivates us to explore a new reformulation for the discretized OCPEC (\ref{equation: discretized OCPEC}) that avoids solving the MPCC.

\section{Proposed reformulation for OCPEC}\label{section: proposed reformulation for OCPEC}

\subsection{Motivation: regularized gap function for VI}
In this section, we propose a \textit{gap-constraint-based reformulation} for the OCPEC.
The core idea of the gap constraint is inspired by the merit-function-based reformulation for the VI.
In this study, we consider a differentiable merit function called the \textit{regularized gap function}, which reformulates the VI as a differentiable \textit{constrained} optimization problem.
The regularized gap function was first introduced independently in \cite{auchmuty1989variational} and \cite{fukushima1992equivalent}.
Its definition is provided below:
\begin{defn}[regularized gap function]
    Let $K \subseteq \mathbb{R}^{n_{\lambda}}$ be a nonempty closed convex set, $F: \mathbb{R}^{n_\lambda} \rightarrow \mathbb{R}^{n_{\lambda}}$ be a continuous function, $c > 0$ be a scalar, and $A \in \mathbb{R}^{n_{\lambda} \times n_{\lambda}}$ be a symmetric positive definite matrix. 
    The \textit{regularized gap function} $\varphi^c: \mathbb{R}^{n_\lambda} \rightarrow \mathbb{R}$ for the VI$(K, F)$ is defined as
    \begin{equation}\label{equation: regularized gap function}
        \varphi^{c}(\lambda) = \max_{\omega \in K} \{F(\lambda)^T (\lambda - \omega) -\frac{c}{2}(\lambda-\omega)^T A (\lambda-\omega) \}.
    \end{equation}
\end{defn}

Given a $\lambda \in \mathbb{R}^{n_\lambda}$, the maximization problem that defines $\varphi^{c}(\lambda)$ is a strictly concave program in the variable $\omega$ and parameterized by $\lambda$. 
Thus, the solution to this program is unique, denoted by $\hat{\omega}$.
Furthermore, it has an explicit expression that is the projection of the stationary point $\lambda - \frac{1}{c}A^{-1}F(\lambda)$ onto the set $K$ under the norm $\| \cdot \|_A$: 
\begin{equation}\label{equation: explicit expression of omega}
    \hat{\omega} = \omega(\lambda) = \Pi_{K, A}(\lambda - \frac{1}{c}A^{-1}F(\lambda)),
\end{equation}
where $\Pi_{K, A}(x)$ is the skewed projector and is defined as the unique solution to a strictly convex program:
\begin{equation}\label{skewed projector}
    \Pi_{K, A}(x) := \arg \min_{y \in K} \{\frac{1}{2}(y - x)^T A (y - x)\}.
\end{equation}
We thereby can write down $\varphi^{c}(\lambda)$ explicitly, provided that $\hat{\omega}$ has been evaluated using (\ref{equation: explicit expression of omega}) based on the given $\lambda$:
\begin{equation}
    \varphi^{c}(\lambda) = F(\lambda)^T (\lambda - \hat{\omega}) -\frac{c}{2}(\lambda-\hat{\omega})^T A (\lambda-\hat{\omega}).
\end{equation}

Here, we focus on two properties of $\varphi^c(\lambda)$: its equivalence to the VI and its differentiability,
as stated in the following theorem (Theorem 10.2.3, \cite{facchinei2003finite}):

\begin{thm}\label{theorem: properties of the regularized gap function}
    The following two statements are valid for a regularized gap function $\varphi^c(\lambda)$ given by (\ref{equation: regularized gap function}):
    \begin{itemize}
        \item (\textit{Equivalence}) $\varphi^c(\lambda) \geq 0, \forall \lambda \in K$. 
        Furthermore, $\varphi^c(\lambda) = 0, \lambda \in K $ if and only if $\lambda \in$ SOL$(K, F)$.
        Hence, $\varphi^c(\lambda)$ is a merit function for VI$(K, F)$ only when $\lambda \in K$, and we can solve VI$(K, F)$ by globally solving a \textit{constrained} optimization problem:
        \begin{equation}\label{equation: constrained optimization problem based on regularized gap function}
             \min_{\lambda \in K} \ \varphi^c(\lambda).
        \end{equation}
        \item (\textit{Differentiability}) If $F$ is continuously differentiable, then so is $\varphi^c(\lambda)$, and $\nabla_{\lambda}\varphi^c(\lambda)$ is explicitly given by
        \begin{equation}\label{equation: explicit expression of regularized gap function}
            \nabla_{\lambda}\varphi^c(\lambda) = F^T(\lambda) + (\lambda - \hat{\omega})^T(\nabla_{\lambda}F (\lambda) - cA).
        \end{equation}
    \end{itemize}
\end{thm}

\begin{rem}\label{remark: stationary properties of regularized gap function}
    Although Theorem \ref{theorem: properties of the regularized gap function} states the necessity of finding the \textit{global minimizer} of (\ref{equation: constrained optimization problem based on regularized gap function}), the equivalence still holds when only the \textit{stationary point} of (\ref{equation: constrained optimization problem based on regularized gap function}) is found, provided that VI$(K, F)$ exhibits certain properties.
    See Theorem 10.2.5 in \cite{facchinei2003finite} for details.
\end{rem}

Based on $\varphi^c(\lambda)$, we propose a new reformulation for the VI, which, in fact, is a variant of the merit-function-based reformulation (\ref{equation: merit function based reformulation for VI}), as stated in the following proposition:
\begin{prop}\label{proposition: gap constraint based reformulation for VI}
    Let $K \subseteq \mathbb{R}^{n_{\lambda}}$ be a nonempty closed convex set given by (\ref{equation: finitely representable VI set K}), $F: \mathbb{R}^{n_\lambda} \rightarrow \mathbb{R}^{n_{\lambda}}$ be a continuous function, and $\varphi^c(\lambda): \mathbb{R}^{n_\lambda} \rightarrow \mathbb{R}$ be a regularized gap function given by (\ref{equation: regularized gap function}).
    We have that $\lambda \in$ SOL$(K, F)$ if and only if $\lambda$ satisfies a set of $n_g + 1$ inequalities:
    \begin{subequations}\label{equation: gap constraint based reformulation for VI}
        \begin{align}
            & g(\lambda) \geq 0, \label{equation: gap constraint based reformulation for VI VI set} \\
            & \varphi^c(\lambda) \leq 0. \label{equation: gap constraint based reformulation for VI gap constraint}
        \end{align}
    \end{subequations}  
    We call (\ref{equation: gap constraint based reformulation for VI}) the \textit {gap-constraint-based reformulation} for the VI$(K, F)$, and (\ref{equation: gap constraint based reformulation for VI gap constraint}) the \textit{gap constraint}.
\end{prop}
\begin{pf}
    The inequality (\ref{equation: gap constraint based reformulation for VI VI set}) implies $\varphi^c(\lambda) \geq 0$; 
    together with inequality (\ref{equation: gap constraint based reformulation for VI gap constraint}), we have $\varphi^c(\lambda) = 0$, and thus the equivalence follows from Theorem \ref{theorem: properties of the regularized gap function}.
\end{pf}

\subsection{Gap-constraint-based reformulation for OCPEC}\label{subsection: gap constraint based reformulation for OCPEC}

Now, we are in a position to state the proposed OCPEC reformulation.
Since the function $F(x, u, \lambda)$ in VI (\ref{equation: discretized OCPEC VI}) also includes variables $x,u$, we introduce an auxiliary variable $\eta = F(x, u, \lambda)$ to reduce the complexity of the regularized gap function.
We then have the following proposition that is inherited from Theorem \ref{theorem: properties of the regularized gap function} and Proposition \ref{proposition: gap constraint based reformulation for VI}:
\begin{prop}\label{proposition: gap constraint based reformulation for OCPEC VI}
    Let $K \subseteq \mathbb{R}^{n_\lambda}$ be a nonempty closed convex set given by (\ref{equation: finitely representable VI set K}), $F: \mathbb{R}^{n_x} \times \mathbb{R}^{n_u} \times \mathbb{R}^{n_\lambda} \rightarrow \mathbb{R}^{n_\lambda}$ be a continuously differentiable function, $c$ be a positive scalar, $A \in \mathbb{R}^{n_{\lambda} \times n_{\lambda}}$ be a symmetric positive definite matrix, and $\eta \in \mathbb{R}^{n_{\lambda}}$ be an auxiliary variable. 
    We define a regularized gap function $\varphi^c: \mathbb{R}^{n_\lambda} \times \mathbb{R}^{n_\lambda} \rightarrow \mathbb{R}$ as
    \begin{equation}\label{equation: regularized gap function OCPEC}
        \varphi^{c}(\lambda, \eta) = \max_{\omega \in K} \{\eta^T (\lambda - \omega) -\frac{c}{2}(\lambda-\omega)^T A (\lambda-\omega) \},
    \end{equation}
    then, the following two statements are valid:
    \begin{itemize}
        \item $\lambda \in$ SOL$(K, F(x, u, \lambda))$ if and only if $(x, u, \lambda, \eta)$ satisfies a set of $n_{\lambda} $ equalities and $n_g + 1$ inequalities:
        \begin{subequations}\label{equation: gap constraint based reformulation for OCPEC VI}
            \begin{align}
                & F(x, u, \lambda) - \eta = 0, \\
                & g(\lambda) \geq 0, \\
                & \varphi^c(\lambda, \eta) \leq 0. \label{equation: gap constraint based reformulation for OCPEC VI gap constraint}   
            \end{align}
        \end{subequations}
        \item $\varphi^{c}(\lambda, \eta)$ is continuously differentiable; moreover, we can explicitly evaluate $\varphi^{c}(\lambda, \eta)$ and its gradient by
        \begin{subequations}\label{equation: evaluating gap constraint for OCPEC}
            \begin{align}
                & \varphi^{c}(\lambda, \eta) = \eta^T (\lambda - \hat{\omega}) -\frac{c}{2}(\lambda-\hat{\omega})^T A (\lambda-\hat{\omega}), \label{equation: evaluating gap constraint for OCPEC, gap function}\\
                & \nabla_{\lambda}\varphi^{c}(\lambda, \eta) = \eta^T -  c(\lambda - \hat{\omega})^T A,\\
                & \nabla_{\eta}\varphi^{c}(\lambda, \eta) = (\lambda - \hat{\omega})^T,
            \end{align}
        \end{subequations}
        with $\hat{\omega}$ being the unique solution to the maximization problem that defines $\varphi^{c}(\lambda, \eta)$:
        \begin{equation}\label{equation: explicit expression of omega OCPEC}
            \hat{\omega} = \omega(\lambda, \eta) = \Pi_{K, A}(\lambda - \frac{1}{c}A^{-1}\eta).
        \end{equation}
    \end{itemize}
\end{prop}

Consequently, by replacing VI (\ref{equation: discretized OCPEC VI}) with its gap-constraint-based reformulation (\ref{equation: gap constraint based reformulation for OCPEC VI}), and further relaxing the gap constraint (\ref{equation: gap constraint based reformulation for OCPEC VI gap constraint}), we have the following \textit{relaxed gap-constraint-based} problem for the discretized OCPEC (\ref{equation: discretized OCPEC}):
\begin{subequations}\label{equation: discretized OCPEC relax gap}
    \begin{align}
    \min_{\{ x_n, u_n,\lambda_n\}^N_{n=1}} & L_T(x_N) + \sum^{N}_{n=1}L(x_n, u_n, \lambda_n), \\
     \text{s.t.}     \quad                 & G(x_n, u_n) \geq 0, \\
                                           & C(x_n, u_n) = 0, \\
                                           & x_{n-1} + \mathcal{F}(x_n, u_n, \lambda_n) = 0, \\
                                           & F(x_n, u_n, \lambda_n) - \eta_n = 0, \label{equation: discretized OCPEC relax gap auxiliary variable}\\
                                           & g(\lambda_n) \geq 0, \\
                                           & s - \varphi^c(\lambda_n, \eta_n) \geq 0, \label{equation: discretized OCPEC relax gap gap constraint}
    \end{align}
\end{subequations}
with $n = 1, \dots ,N$, $s \geq 0$ the relaxation parameter, and $\eta_n \in \mathbb{R}^{n_\lambda}$ the auxiliary variable for $F(x_n, u_n, \lambda_n)$.

Even introducing $n$ auxiliary variables $\eta_n$, the new relaxed problem (\ref{equation: discretized OCPEC relax gap}) still possesses a more concise constraint system: the number of constraints related to equilibrium constraints is only $N(n_{\lambda} + n_g + 1)$ (i.e., (\ref{equation: discretized OCPEC relax gap auxiliary variable})--(\ref{equation: discretized OCPEC relax gap gap constraint})), 
whereas in the relaxed problem (\ref{equation: discretized OCPEC relax KKT}), this number is $N(n_{\lambda} + 3n_g)$ (i.e., (\ref{equation: discretized OCPEC relax KKT VI stationary})--(\ref{equation: discretized OCPEC relax KKT VI complementary})).
Furthermore, this new relaxed problem (\ref{equation: discretized OCPEC relax gap}) exhibits two additional favorable properties: 
First, it is a differentiable problem benefiting from the continuous differentiability of $\varphi^c(\lambda, \eta)$. 
Second, its feasible set is equivalent to that of the original problem (\ref{equation: discretized OCPEC}) when $s= 0$, and possesses a feasible interior when $s > 0$, as shown in Section \ref{section: geometric interpretation}.
Thus, solutions to the discretized OCPEC (\ref{equation: discretized OCPEC}) can also be obtained by solving a sequence of relaxed problems (\ref{equation: discretized OCPEC relax gap}) with $s \rightarrow 0$.

\begin{rem}\label{remark: constraint qualification discussion}
    It is crucial to investigate whether the new relaxed problem (\ref{equation: discretized OCPEC relax gap}) satisfies the constraint qualification when $s = 0$.
    Unfortunately, the answer is negative.
    In Theorem \ref{theorem: properties of the regularized gap function}, the zeros of $\varphi^c(\lambda)$ within the set $K$ are also the global solutions to the constrained optimization problem (\ref{equation: constrained optimization problem based on regularized gap function}).
    Consequently, for any feasible point that satisfies (\ref{equation: gap constraint based reformulation for VI}), the gap constraint $\varphi^c(\lambda) \leq 0$ must be active, and the gradient of $\varphi^c(\lambda)$ is either zero or linearly dependent with the gradient of the activated $g(\lambda) \geq 0$, which essentially violates the constraint qualification. 
    A similar discussion about the bilevel optimization can be found in \cite{ouattara2018duality}, where the constraint qualification is interpreted as stating the constraints without the optima of an embedded optimization problem.
\end{rem}

\section{Proposed solution method for OCPEC}\label{section: proposed solution method for OCPEC}

\subsection{Efficient evaluation of gap constraints}\label{subsection: efficient evaluation of gap constraints}
While the gap-constraint-based reformulation exhibits numerous favorable theoretical properties, the evaluation of gap functions is far from straightforward, and one could argue that this is indeed the main drawback of this reformulation.
As demonstrated in (\ref{equation: evaluating gap constraint for OCPEC}), each computation of the gap function requires solving one constrained optimization problem, which seriously hinders the practical application.
Nonetheless, evaluating the skewed projector (\ref{equation: explicit expression of omega OCPEC}) can be quite efficient in some common cases.
For example, if we specify matrix $A$ as an identity matrix, and set $K$ exhibits a box-constraint structure:
\begin{equation}\label{equation: box constraint VI set K}
    K := \{\lambda \in \mathbb{R}^{n_\lambda} | b_l \leq \lambda \leq b_u \},
\end{equation}
with $b_l \in \{\mathbb{R} \cup \{-\infty\} \}^{n_{\lambda}}$, $b_u \in \{\mathbb{R} \cup \{+\infty\} \}^{n_{\lambda}}$, and $b_l < b_u$,
then, the projection operator in (\ref{equation: explicit expression of omega OCPEC}) requires only one call to the max operator and one call to the min operator:
\begin{equation}\label{equation: evaluate omega BVI case}
    \hat{\omega} = \omega(\lambda, \eta) = \min(\max(b_l, \lambda - \frac{1}{c}\eta), b_u),
\end{equation}
which involves a very low computational cost.

Even in the absence of the box-constraint structure in set $K$, the skewed projector (\ref{equation: explicit expression of omega OCPEC}) can still be computationally tractable by exploiting the OCPEC structure.
First, it is noteworthy that we can evaluate $\varphi^c(\lambda_n, \eta_n)$ in a parallel fashion with up to $N$ cores, as $\varphi^c(\lambda_n, \eta_n)$ only depends on the variables of the stage $n$.
Second, even if we evaluate $\varphi^c(\lambda_n, \eta_n)$ in a serial fashion, note that the maximization problem that defines $\varphi^c(\lambda_n, \eta_n)$ is only parameterized by $\lambda_n, \eta_n$, and if the parameters of adjacent problems do not change significantly, we can expect that the optimal active sets of the adjacent problems exhibit slight differences, or even remain unchanged.
This enables us to efficiently evaluate the skewed projectors (\ref{equation: explicit expression of omega OCPEC}) by certain active set warm-start techniques. 
For example, if set $K$ is polyhedral, then the skewed projector (\ref{equation: explicit expression of omega OCPEC}) is reduced to the solution to a strictly concave quadratic program (QP):
\begin{equation}\label{equation: evaluate omega QP case}
    \hat{\omega} = \omega(\lambda, \eta) = \arg \max_{\omega \in K}\{-\frac{c}{2}\omega^T A \omega - (\eta^T - c\lambda^TA)\omega \}.
\end{equation}
In this scenario, a suitable choice is the solver qpOASES \citep{Ferreau2014}, an efficient small-scale QP solver using an \textit{online active set strategy} \citep{Ferreau2008}.

\subsection{Successive gap constraint linearization method}
Although we cannot directly use off-the-shelf NLP solvers to solve the relaxed problem (\ref{equation: discretized OCPEC relax gap}) owing to the particular nature of evaluating $\varphi^c(\lambda_n, \eta_n)$, we can still utilize the NLP algorithm framework to solve the problem (\ref{equation: discretized OCPEC relax gap}).
Here, we adopt the \textit{sequential quadratic programming} (SQP) framework, which is a Newton-type method and facilitates the use of state-of-the-art QP solvers.

We first make some modifications to the relaxed problem (\ref{equation: discretized OCPEC relax gap}).
We introduce scalar auxiliary variables $v_n$ for the regularized gap functions $v_n = \varphi^c(\lambda_n, \eta_n)$ to reduce the nonlinearity of inequality constraints and facilitate easier initialization. 
We also penalize $\varphi^c(\lambda_n, \eta_n)$ \textit{indirectly} by adding a quadratic penalty term for the variable $v_n$ with a penalty parameter $\mu > 0$, as stated in (\ref{equation: compact NLP cost function}). 
By collecting 
all decision variables into a vector $\boldsymbol{z} = [z^T_1, \cdots z^T_n, \cdots z^T_N]^T$ with $z_n = [x^T_n, u^T_n, \lambda^T_n, \eta^T_n, v_n]^T$, 
all equality constraints into a vector $\boldsymbol{h} = [h^T_1, \cdots h^T_n, \cdots h^T_N]^T$ with $h_n = [x^T_{n-1} + \mathcal{F}^T_n, C^T_n, F^T_n - \eta^T_n, \varphi^c_n - v_n]^T$, 
and all inequality constraints into a vector $\boldsymbol{c} = [c^T_1, \cdots c^T_n, \cdots c^T_N]^T$ with $c_n = [G^T_n, g^T_n, s-v_n]^T$, 
we can rewrite problem (\ref{equation: discretized OCPEC relax gap}) as a general NLP:
\begin{subequations}\label{equation: compact NLP}
    \begin{align}
        \min_{\boldsymbol{z}} \ & \boldsymbol{J}(\boldsymbol{z}, \mu), \\
        \text{s.t.}     \quad   & \boldsymbol{h}(\boldsymbol{z}) = 0, \\
                                & \boldsymbol{c}(\boldsymbol{z}, s) \geq 0,
    \end{align}
\end{subequations}
with the cost function being
\begin{equation}\label{equation: compact NLP cost function}
    \boldsymbol{J}(\boldsymbol{z}, \mu) = L_T(x_N) + \sum^{N}_{n=1}\{L(x_n, u_n, \lambda_n) + \mu \Delta t v^2_n\}.
\end{equation}
In the relaxed problem (\ref{equation: compact NLP}), the regularized gap function $\varphi^c(\lambda_n, \eta_n)$ is simultaneously relaxed and penalized.
On one hand, the constraint $s-v_n \geq 0$ serves as a \textit{hard constraint} that provides an upper bound on $\varphi^c(\lambda_n, \eta_n)$. 
On the other hand, the quadratic penalty term $\mu \Delta t v^2_n$ serves as a \textit{soft constraint} that attempts to drive $\varphi^c(\lambda_n, \eta_n)$ to zero.

We now commence the introduction of the proposed SQP-type method, called the \textit{successive gap constraint linearization} (SGCL) method, to solve (\ref{equation: compact NLP}) with a given $s > 0$.
Algorithm \ref{algorithm: SGCL} summarizes the SGCL method.
The primary difference between our SGCL method and classical SQP methods lies in step 1, namely, for the current iterate $\boldsymbol{z}^k$, it is necessary to \textit{first} evaluate skewed projectors $\{\hat{\omega}^k_n\}^N_{n=1}$ in a stage-wise manner (Lines 4--6) for the computation of the regularized gap function and its gradient (Line 8).

\begin{algorithm}[!tbp]
\caption{SGCL method to solve the relaxed problem (\ref{equation: compact NLP}) with a given $s > 0$}
\label{algorithm: SGCL}
\begin{algorithmic}
    \State \textbf{Input: } $s$, $\boldsymbol{z}^0$
    \State \textbf{Output: } $\boldsymbol{z}^*$ 
\end{algorithmic}

\begin{algorithmic}[1]
    \State \textbf{Initialization: } $\boldsymbol{z}^k \gets \boldsymbol{z}^0, \boldsymbol{\gamma}^k_h, \boldsymbol{\gamma}^k_c \gets 0$, select $\mu > 0$.
    \For{$k = 1 \text{ to } k_{max}$}
    \State \textbf{Step 1. } \text{Evaluate functions and derivatives: }
        \For{$n = 1 \text{ to } N $} \text{in serial or parallel}
            \State $\hat{\omega}^{k}_n \gets \text{extract } \lambda^{k}_n, \eta^{k}_n \text{ from } \boldsymbol{z}^k  \text{ and use } (\ref{equation: explicit expression of omega OCPEC})$.
        \EndFor    
    \State $\boldsymbol{J}^{k}, \boldsymbol{c}^{k}, \nabla_{\boldsymbol{z}} \boldsymbol{J}^{k}, \nabla_{\boldsymbol{z}}\boldsymbol{c}^{k}, \boldsymbol{H}^{k} \gets \boldsymbol{z}^{k}, \mu, s$.
    \State $\boldsymbol{h}^{k}, \nabla_{\boldsymbol{z}} \boldsymbol{h}^{k} \gets \boldsymbol{z}^{k}, \{\hat{\omega}^{k}_n\}^N_{n=1}$.
    \State \textbf{Step 2. } \text{Evaluate search direction and multiplier: }
    \State $\Delta \boldsymbol{z}, \hat{\boldsymbol{\gamma}}^{k+1}_h, \hat{\boldsymbol{\gamma}}^{k+1}_c\gets  \text{solve a sparse convex QP}$ (\ref{equation: qp for search direction}).
    \State \textbf{Step 3. } \text{Check termination condition: }
    \If{ one of (\ref{equation: terminal condition KKT error}) (\ref{equation: terminal condition norm of search direction}) (\ref{equation: terminal condition decrease in gap function}) is satisfied}
        \State $\boldsymbol{z}^* \gets \boldsymbol{z}^k$ and break iteration routine.
    \EndIf
    \State \textbf{Step 4. } \text{Filter line search globalization: }
    \State $\boldsymbol{z}^{k+1}, \alpha \gets \boldsymbol{z}^k, \Delta \boldsymbol{z}$.
    \State \textbf{Step 5. } \text{Update Lagrangian multipliers: }
    \State $\boldsymbol{\gamma}^{k+1}_h \gets \boldsymbol{\gamma}^k_h + \alpha (\hat{\boldsymbol{\gamma}}^{k+1}_h - \boldsymbol{\gamma}^k_h)$.
    \State $\boldsymbol{\gamma}^{k+1}_c \gets \boldsymbol{\gamma}^k_c + \alpha (\hat{\boldsymbol{\gamma}}^{k+1}_c - \boldsymbol{\gamma}^k_c)$.
    \EndFor
\end{algorithmic}
\end{algorithm}

Let the Lagrangian of the relaxed problem (\ref{equation: compact NLP}) be
\begin{equation}
    \mathcal{L}(\boldsymbol{z}, \boldsymbol{\gamma}_h, \boldsymbol{\gamma}_c, \mu, s) = J(\boldsymbol{z}, \mu) + \boldsymbol{\gamma}_h^T\boldsymbol{h}(\boldsymbol{z}) - \boldsymbol{\gamma}^T_c \boldsymbol{c}(\boldsymbol{z}, s),
\end{equation}
with $\boldsymbol{\gamma}_h \in \mathbb{R}^{n_h}, \boldsymbol{\gamma}_c \in \mathbb{R}^{n_c}$ being the Lagrangian multipliers for constraints $\boldsymbol{h}$ and $\boldsymbol{c}$, respectively.
Let $\boldsymbol{H}$ be the Hessian of the Lagrangian $\mathcal{L}$.
Then, we can obtain a new search direction $\Delta \boldsymbol{z}$ and an estimate of the new multipliers $\hat{\boldsymbol{\gamma}}^{k+1}_h$, $\hat{\boldsymbol{\gamma}}^{k+1}_c$ by solving the following structured sparse QP:
\begin{subequations}\label{equation: qp for search direction}
    \begin{align}
        \min_{\Delta\boldsymbol{z}} \ & \frac{1}{2} \Delta\boldsymbol{z}^T \boldsymbol{H}^k \Delta\boldsymbol{z} + \nabla_{\boldsymbol{z}}\boldsymbol{J}^k \Delta\boldsymbol{z}, \\
        \text{s.t.}     \quad   & \boldsymbol{h}^k + \nabla_{\boldsymbol{z}} \boldsymbol{h}^k \Delta\boldsymbol{z} = 0, \\
                                & \boldsymbol{c}^k + \nabla_{\boldsymbol{z}} \boldsymbol{c}^k \Delta\boldsymbol{z} \geq 0.        
    \end{align}
\end{subequations}
In practice, the cost terms $L_T, L$ in (\ref{equation: compact NLP cost function}) are often formulated as quadratic terms.
In this scenario, we can ensure the convexity of QP (\ref{equation: qp for search direction}) by the Gauss--Newton approximation that exploits only the convexity of the cost $\boldsymbol{J}$:
\begin{equation}
    \boldsymbol{H}^k \approx \nabla_{\boldsymbol{z}\boldsymbol{z}} \boldsymbol{J}^k \succeq 0.
\end{equation}
Exploiting convexity for other problem structures can be found in an excellent survey \citep{messerer2021survey}.

After obtaining the search direction $\Delta \boldsymbol{z}$, we decide whether to terminate the iteration routine by examining certain quantities.
One of these quantities is the optimality errors, including the primal residual $E^k_p$, dual residual $E^k_d$, and complementary residual $E^k_c$.
Their definitions are below:
\begin{subequations}\label{equation: optimality error}
    \begin{align}
        & E^k_p = \max (\| \boldsymbol{h}^k\|_{\infty}, \| \min(0, \boldsymbol{c}^k) \|_{\infty}), \\
        & E^k_d = \frac{1}{\kappa_d}\max (\|\nabla_{\boldsymbol{z}} \mathcal{L}^k\|_{\infty}, \| \min(0, \boldsymbol{\gamma}^k_c) \|_{\infty}), \\
        & E^k_c = \frac{1}{\kappa_c} \| \boldsymbol{c}^k \odot \boldsymbol{\gamma}^k_c  \|_{\infty},
    \end{align}
\end{subequations}
with $\kappa_d$, $\kappa_c \geq 1$ the scaling parameters defined in \cite{wachter2006implementation}.
The iteration routine will be terminated if any of the following three conditions are satisfied:
\begin{itemize}
    \item Condition based on a composite (KKT) residual:
    \begin{equation}\label{equation: terminal condition KKT error}
        \max (E^k_p, E^k_d, E^k_c) \leq \epsilon_{kkt};
    \end{equation}
    \item Condition based on the search direction:
    \begin{equation}\label{equation: terminal condition norm of search direction}
         \| \Delta \boldsymbol{z} \|_{\infty} \leq \epsilon_{sd};
    \end{equation}
    \item Condition based on the individual residuals:
    \begin{equation}\label{equation: terminal condition decrease in gap function}
        E^k_p \leq \epsilon_p, \quad E^k_d \leq \epsilon_d, \quad E^k_c \leq \epsilon_c, 
    \end{equation}
\end{itemize}
where $\epsilon_{kkt}, \epsilon_{sd}, \epsilon_p, \epsilon_d,  \epsilon_c > 0$ are the given tolerances. 

In case the terminal conditions are not satisfied, we use the filter line search globalization method \citep{wachter2006implementation} to find a new iterate $\boldsymbol{z}^{k+1} = \boldsymbol{z}^{k} + \alpha \Delta \boldsymbol{z}$ with stepsize $\alpha$. 
$\boldsymbol{z}^{k+1}$ achieves decrease either in the cost $\boldsymbol{J}(\boldsymbol{z}, \mu)$ or the constraint violation $\boldsymbol{M}(\boldsymbol{z}, s)$, where $\boldsymbol{M}(\boldsymbol{z}, s)$ is defined as the $\ell_1$ norm of the constraint scaled by $\Delta t$:
\begin{equation}\label{equation: constraint violation}
    \boldsymbol{M}(\boldsymbol{z}, s) = \Delta t (\| \boldsymbol{h} \|_1 + \| \min(0, \boldsymbol{c})\|_1).
\end{equation}
We sometimes re-initialize the filter as a heuristic \citep{wachter2006implementation} for better performance.

As discussed in Remark \ref{remark: constraint qualification discussion} about the gradient of $\varphi^{c}(\lambda, \eta)$, the SGCL method will also encounter numerical difficulties when $s$ is close to zero.
These numerical difficulties mainly arise from three aspects:
First, QP (\ref{equation: qp for search direction}) numerically violates constraint qualifications.
Second, the search direction $\Delta \boldsymbol{z}$ is insufficient to reduce the constraint violation, as the gradient of $\varphi^{c}(\lambda, \eta)$ may approach zero.
Third, the dual residual $E^k_d$ fails to satisfy a composite tolerance $\epsilon_{kkt}$, as the elements in $\boldsymbol{\gamma}^k_h, \boldsymbol{\gamma}^k_c$ associated with $\varphi^{c}(\lambda, \eta), g(\lambda)$ may become extremely large.
We mitigate these numerical difficulties using the following two approaches:
We solve QP (\ref{equation: qp for search direction}) using a robust and efficient solver FBstab \citep{liao2020fbstab}, which needs neither constraint qualifications nor a feasible interior.
We carefully select $\epsilon_p, \epsilon_d,  \epsilon_c$ in (\ref{equation: terminal condition decrease in gap function}) to ensure that the algorithm can still be terminated promptly on an optimal solution.

\begin{rem}\label{remark: convex approximation}
    Our SGCL method is named based on the fact that gap functions appear as the problem constraint rather than being a part of the cost function, and we only make use of its function value and gradient.  
    The difficulty in approximating $\boldsymbol{H}^k$ is the major reason why we do not directly penalize $\varphi^c(\lambda, \eta)$ but instead indirectly penalize it, as $\varphi^c(\lambda, \eta)$ is generally nonconvex and only once continuously differentiable. 
    Currently, we are exploring methods to efficiently construct a convex approximation for $\varphi^c(\lambda, \eta)$ and hope to exploit its second-order information.
\end{rem}

\section{Geometric interpretation}\label{section: geometric interpretation}

\begin{figure}[!tbp]
    \centering
    \includegraphics[width=0.8\linewidth]{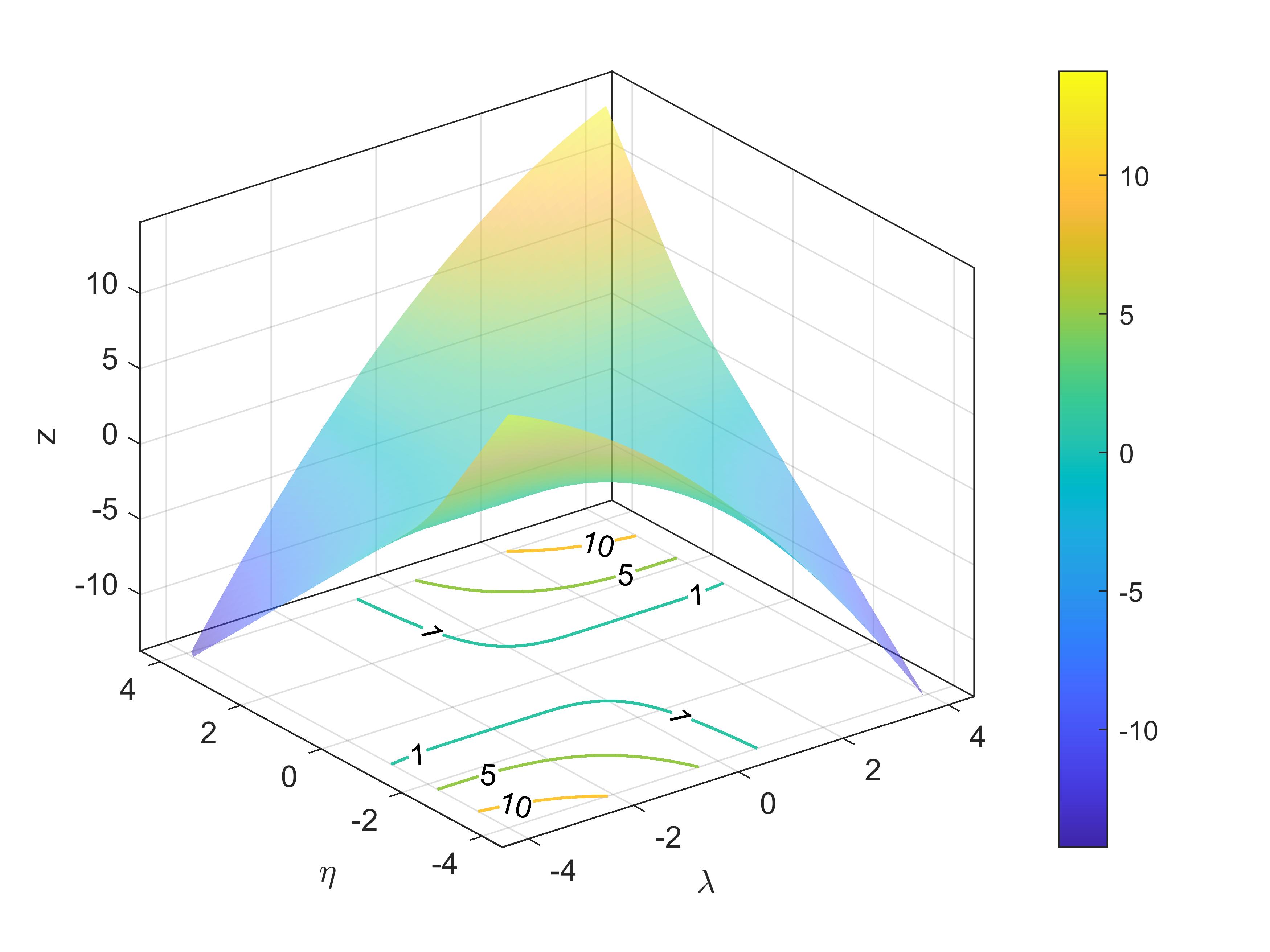}
    \caption{Contour of $\varphi^c(\lambda, \eta)$ with $c = 0.5, b_l = -1, b_u = 1$.}
    \label{fig: regularized gap function contour}
\end{figure}
\begin{figure}[!tbp]
    \centering
    \includegraphics[width=0.8\linewidth]{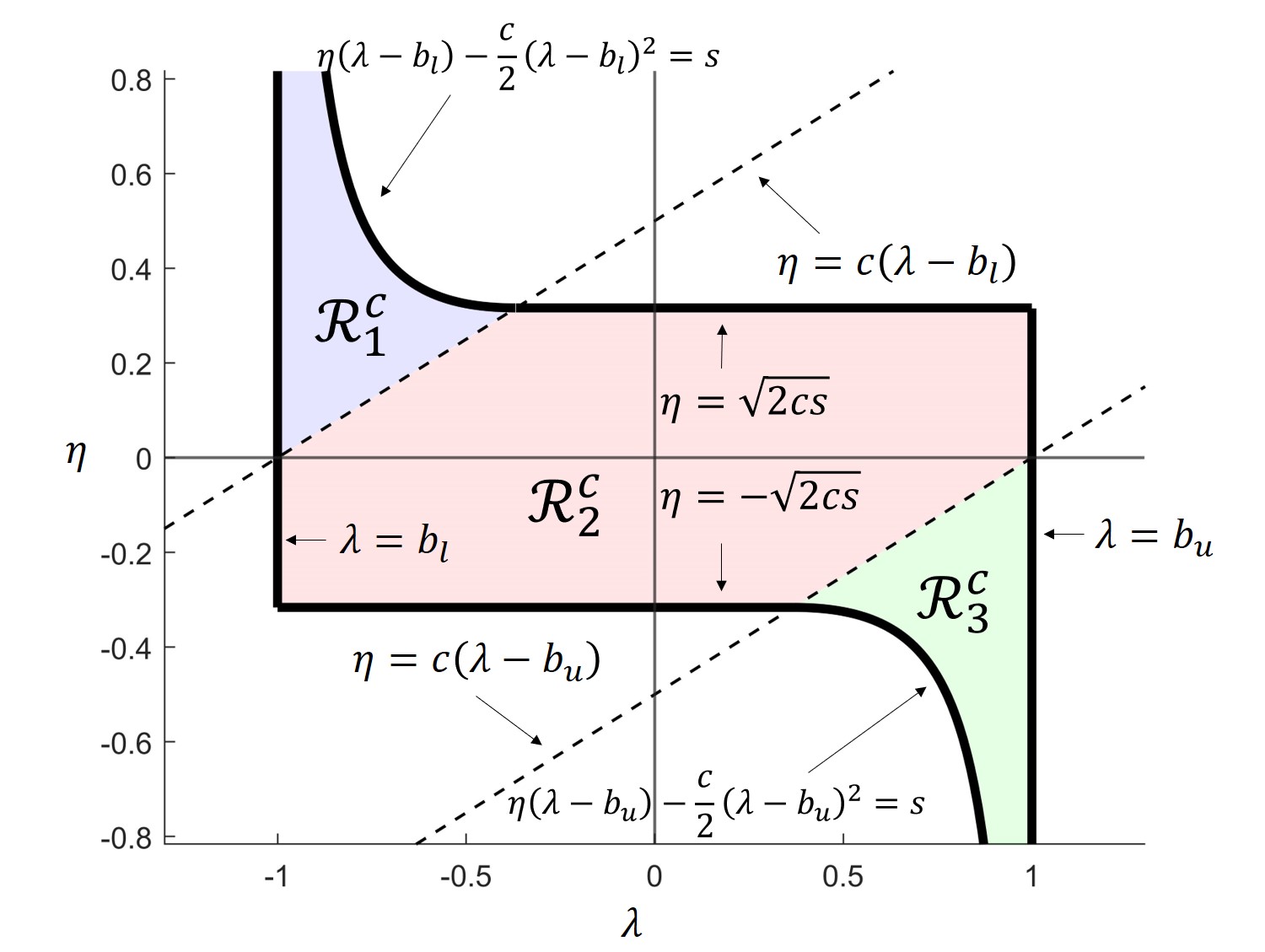}
    \caption{Relaxed feasible set $\mathcal{R}^c  = \mathcal{R}^c_1 \bigcup \mathcal{R}^c_2 \bigcup \mathcal{R}^c_3$, with $c = 0.5$, $s = 0.1$, $b_l = -1$, and $b_u = 1$. }
    \label{fig: regularized gap function region}
\end{figure}

In this section, we provide an intuitive geometric interpretation of the gap constraint using a simple mathematical programming with equilibrium constraint (MPEC):
\begin{subequations}\label{equation: simple MPEC}
    \begin{align}
        \min_{\lambda, \eta} \ & J(\lambda, \eta), \\
        \text{s.t. } &  \lambda \in \text{SOL}(K, \eta), \label{equation: simple MPEC constraint}
    \end{align}
\end{subequations}
where $\lambda, \eta \in \mathbb{R}$ are scalar variables to be optimized, set $K$ exhibits a box-constrained structure in the form of (\ref{equation: box constraint VI set K}), and $J: \mathbb{R} \times \mathbb{R} \rightarrow \mathbb{R}$ is the smooth cost function.
The feasible set formed by (\ref{equation: simple MPEC constraint}) includes three pieces: the nonnegative part of axis $\lambda = b_l$, the nonpositive part of axis $\lambda = b_u$, and the segment of axis $\eta = 0$ between $\lambda = b_l$ and $\lambda = b_u$.
By specifying $A = I_{1 \times 1}$ and substituting (\ref{equation: evaluate omega BVI case}) into (\ref{equation: evaluating gap constraint for OCPEC, gap function}), we have a relaxed gap-constraint-based problem for MPEC (\ref{equation: simple MPEC}) with a relaxation parameter $s \geq 0$:
\begin{subequations}\label{equation: relax gap constraint based problem for simple MPEC}
    \begin{align}
        \min_{\lambda, \eta} \ & J(\lambda, \eta), \\
        \text{s.t. } &  b_l \leq \lambda \leq b_u, \label{equation: relax gap constraint based problem for simple MPEC set K} \\
                     &  s - \varphi^c(\lambda, \eta) \geq 0, \label{equation: relax gap constraint based problem for simple MPEC gap constraint} 
    \end{align}
\end{subequations}
where $\varphi^c(\lambda, \eta)$ has an explicit expression:
\begin{equation}\label{equation: relax gap constraint based problem for simple MPEC gap constraint explicit}
    \varphi^c(\lambda, \eta) =  
        \begin{cases}
            \eta(\lambda - b_l) - \frac{c}{2}(\lambda - b_l)^2, & \lambda - \frac{1}{c}\eta \leq b_l,\\
            \frac{1}{2c}\eta^2,                                 & b_l \leq \lambda - \frac{1}{c}\eta \leq b_u, \\
            \eta(\lambda - b_u) - \frac{c}{2}(\lambda - b_u)^2, & \lambda - \frac{1}{c}\eta \geq b_u,
        \end{cases}     
\end{equation}
and its contour is shown in Fig. \ref{fig: regularized gap function contour}.
Therefore, the relaxed feasible set $\mathcal{R}^c$ formed by (\ref{equation: relax gap constraint based problem for simple MPEC set K}) and (\ref{equation: relax gap constraint based problem for simple MPEC gap constraint}) with $s \geq 0$ is divided into three subsets $\mathcal{R}^c_1$, $\mathcal{R}^c_2$, and $\mathcal{R}^c_3$ along axes $\eta = c(\lambda - b_l)$ and $\eta = c(\lambda - b_u)$ such that $\mathcal{R}^c = \mathcal{R}^c_1 \bigcup \mathcal{R}^c_2 \bigcup \mathcal{R}^c_3$.
Here, $\mathcal{R}^c_1$, $\mathcal{R}^c_2$, $\mathcal{R}^c_3$ are respectively represented by
\begin{subequations}
    \begin{align}
        \begin{split}
            \mathcal{R}^c_1 = \{ (\lambda, \eta)| &  b_l \leq \lambda \leq b_u, \eta \geq c(\lambda - b_l), \\
            & \eta(\lambda - b_l) - \frac{c}{2}(\lambda - b_l)^2 \leq s\},
        \end{split}
        \\
        \begin{split}
            \mathcal{R}^c_2 = \{(\lambda, \eta)| & b_l \leq \lambda \leq b_u, | \eta | \leq \sqrt{2cs}, \\
            & c(\lambda - b_u) \leq \eta \leq c(\lambda - b_l)\},
        \end{split}
        \\
        \begin{split}
            \mathcal{R}^c_3 = \{(\lambda, \eta)| & b_l \leq \lambda \leq b_u, \eta \leq c(\lambda - b_u), \\
            & \eta(\lambda - b_u) - \frac{c}{2}(\lambda - b_u)^2 \leq s\},
        \end{split} 
    \end{align}
\end{subequations}
and shown in Fig. \ref{fig: regularized gap function region}, where the blue, red, and green regions are $\mathcal{R}^c_1$, $\mathcal{R}^c_2$, and $\mathcal{R}^c_3$, respectively.
Clearly, $\mathcal{R}^c$ exhibits two prominent features:
First, it always contains the feasible region of the original MPEC (\ref{equation: simple MPEC}), regardless of the choice of $s \geq 0$.
Second, it is a \textit{connected} set, and thus we can expect that the relaxed problem (\ref{equation: relax gap constraint based problem for simple MPEC}) with a given $s > 0$ can be easily solved by the NLP solving algorithm.

\section{Numerical experiment}\label{section: numerical experiment}

\subsection{Benchmark problem and implementation details}
We consider OCPEC (\ref{equation: continuous OCPEC}) with a quadratic cost function:
\begin{subequations}
    \begin{align*}
        L_{T} = \ & \| x(T) - x_{e} \|^2_{Q_T},
        \\
        L_{S} = \ & \| x(t) - x_{ref}(t) \|^2_{Q_x} + \|u(t)\|^2_{Q_{u}} + \|\lambda(t) \|^2_{Q_{\lambda}},
    \end{align*}
\end{subequations}
where $Q_T, Q_x, Q_{u},Q_{\lambda}$ are weighting matrices, $x_{e}$ is the terminal state, and $x_{ref}(t)$ is the reference trajectory.
We consider an affine DVI generalized from an example of the linear complementary system in \cite{Vieira2020}:
\begin{subequations}\label{equation: affine DVI}
    \begin{align*}
        & \Dot{x} =  \begin{bmatrix} 1 & -3 \\ -8 & 10\end{bmatrix} x  + \begin{bmatrix} 4 \\ 8\end{bmatrix} u + \begin{bmatrix} -3 \\ -1\end{bmatrix} \lambda, \\
        & F = \begin{bmatrix} 1 & -3\end{bmatrix} x + 3 u + 5 \lambda, \\
        & K = \{\lambda \ |  -1 \leq \lambda \leq 1 \}.
    \end{align*}
\end{subequations}
The goal is to drive the state $x$ from $x_{0} = [-\frac{1}{2}, -1]^T$ to $x_{e} = [0, 0]^T$ using the control $u$.
We specify $x_{ref}(t) \equiv x_{e}$, $T = 1 s$, $N = 100$, and enforce bounds on the variables $-2 \leq x_1, x_2 \leq 2 $ and $ -2 \leq u \leq 2$.

The SGCL method is implemented in MATLAB R2023b using the CasADi symbolic framework \citep{Andersson2019}.
All the experiments were performed on a laptop with a 1.80 GHz Intel Core i7-8550U.
We solve the sparse QP (\ref{equation: qp for search direction}) by the sparse implementation of the FBstab with its default settings.
We specify $c = 1$, $A = I_{n_{\lambda} \times n_{\lambda}}$, and tolerance $\epsilon_{kkt} = 10^{-5}$, $\epsilon_d = 10^{-4}$, $\epsilon_{sd}, \epsilon_p, \epsilon_c = 10^{-8}$.

\subsection{Comparison between various SGCL implementations}
Since the set $K$ in the affine DVI example can be viewed as having both a box-constraint structure and a polyhedral structure, the SGCL implementation is categorized into two types based on the evaluation of the skewed projector (\ref{equation: explicit expression of omega OCPEC}).
One utilizes the explicit expression (\ref{equation: evaluate omega BVI case}), denoted by SGCL-box. 
The other solves the concave QP (\ref{equation: evaluate omega QP case}) using qpOASES in a serial fashion, denoted by SGCL-polyh.
Both implementations use the CasADi construct \textit{map} to speed up the for-loop routine (Lines 4--6 in Algorithm \ref{algorithm: SGCL}).

\begin{figure}[!tbp]
    \centering
    \includegraphics[width=0.9\linewidth]{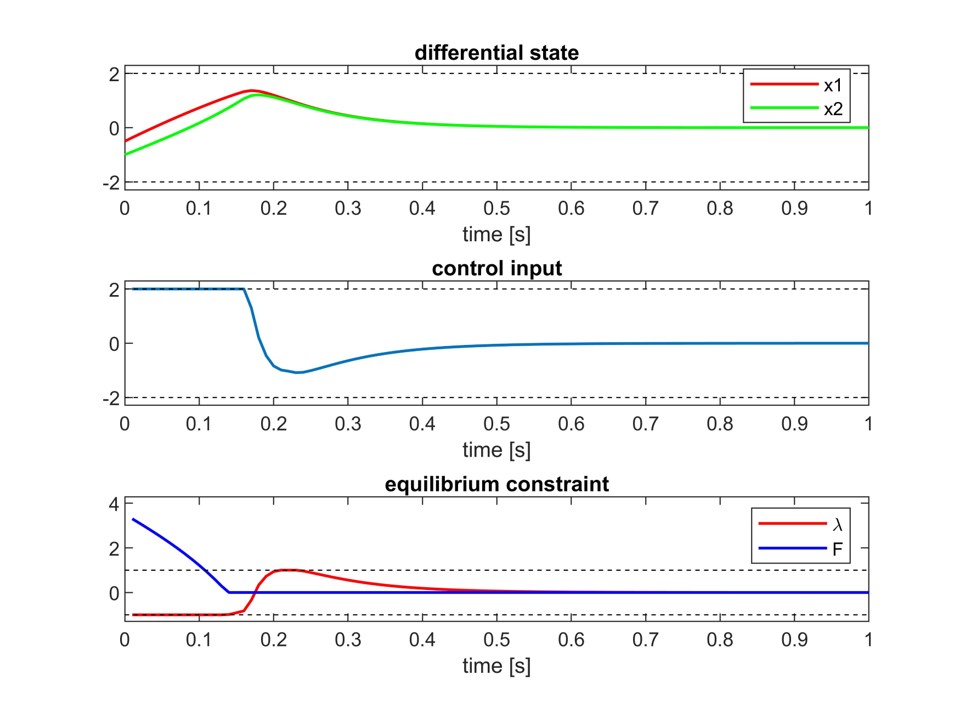}
    \caption{Solution trajectory of the affine DVI example.}
    \label{fig: affine DVI solution trajectory}
\end{figure}

Fig. \ref{fig: affine DVI solution trajectory} illustrates a solution trajectory obtained by solving (\ref{equation: compact NLP}) with a given $s = 10^{-6}$.
Both SGCL-box and SGCL-polyh implementations can obtain this solution after 14 iterations, and the only distinction is the time of evaluating $\hat{\omega}$. 
The majority of the computation time for both implementations is dedicated to solving the sparse QP (\ref{equation: qp for search direction}) as indicated in Table \ref{table: Computation time of the SGCL implementation}.
The time spent on evaluating $\hat{\omega}$ in SGCL-box is negligible, accounting for $\frac{6.11}{214.09} \approx 2.9 \% $ of the total time, which is very promising.
The SGCL-polyh, even with the requirement to solve a sequence of QP for evaluating $\hat{\omega}$, is still efficient, and the time spent on evaluating $\hat{\omega}$ accounts for $\frac{26.56}{235.13} \approx 11.3 \% $ of the total time. 
This validates our discussion in Subsection \ref{subsection: efficient evaluation of gap constraints}, where it was stated that evaluating gap constraints in the OCPEC is not a significant challenge.

\begin{table}[!tbp]
    \caption{Average computation time [ms] of various SGCL implementations, evaluated by running 100 repeated tests initialized with $\boldsymbol{z}^0$ a unit vector}
    \label{table: Computation time of the SGCL implementation}
    \centering
    \begin{tabularx}{\columnwidth} { 
    >{\centering\arraybackslash}X
    >{\centering\arraybackslash}X 
    >{\centering\arraybackslash}X}
    \hline
                                                                            &  SGCL-box    &   SGCL-polyh   \\
    \hline 
    evaluate $\hat{\omega}$ (\ref{equation: explicit expression of omega OCPEC})  &  \textbf{6.11}           &  \textbf{26.56}              \\
    evaluate derivative                                                     &     16.42         &      16.45          \\
    solve QP (\ref{equation: qp for search direction})                      &     174.32         &     174.91           \\
    filter line search                                                      &     15.68         &      15.66          \\
    else                                                                    &     1.55         &      1.56          \\
    total                                                                   &     \textbf{214.09}         &     \textbf{235.13}           \\
    \hline
    \end{tabularx}
\end{table}

\subsection{Comparison of the SGCL and MPCC methods}
We also compare the SGCL method (SGCL-box) with the MPCC method, that is, solving the relaxed problem (\ref{equation: discretized OCPEC relax KKT}) using the NLP solver IPOPT \citep{wachter2006implementation} with its default settings (except the KKT tolerance being $10^{-5}$).
Both methods employ the continuation method that tracks a decreasing parameter sequence $\{ s^j\}_{j = 0}^{J}$, with $s^{j+1} = \max\{s^J, \min \{\kappa_t s^j, (s^j)^{\kappa_e}\} \}$ and $s^0 = 10^{-1}, \kappa_t = 0.8, \kappa_e = 1.5$.
We compare the solutions found by these methods under various $s^J$ in terms of the optimal cost, equilibrium constraint violation, and computation time.
For a given $s^J$, each method is initialized with 20 random $\boldsymbol{z}^0$ that enables the method to find an optimal solution.
The violation of equilibrium constraints is measured using the \textit{natural residual} $\Phi = \lambda - \Pi_{K}(\lambda - F)$ with the Euclidean projector $\Pi_{K}$, and $\Phi = 0$ if and only if $\lambda \in \text{SOL}(K, F)$.

\begin{figure}[!tbp]
    \centering
    \includegraphics[width=0.8\linewidth]{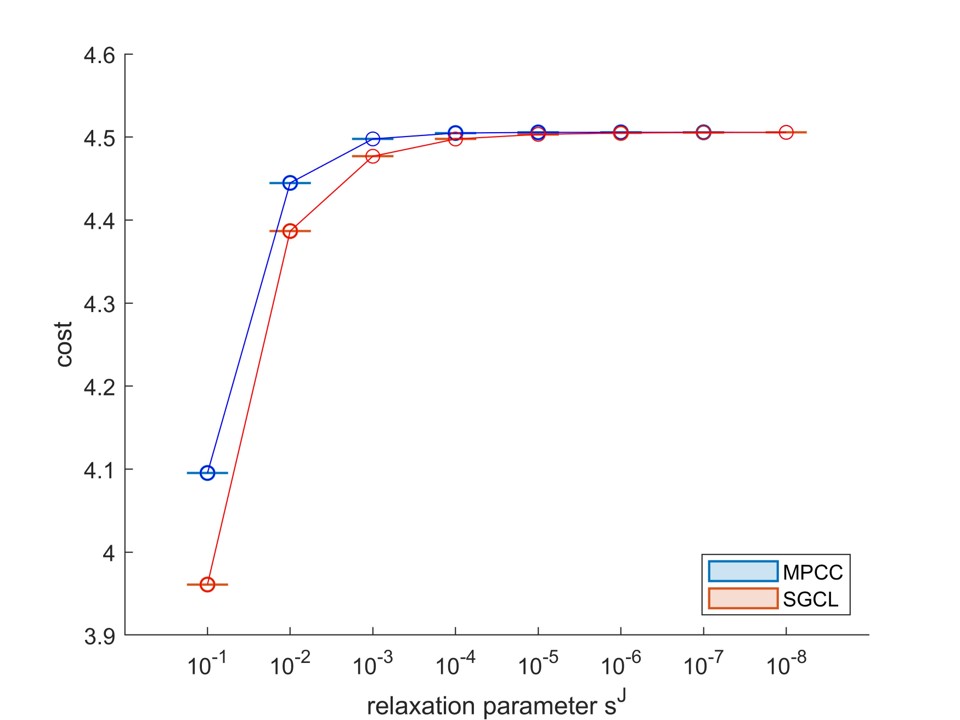}
    \caption{Box plot for relaxation parameter $s^J$ vs. cost.}
    \label{fig: param vs cost}
\end{figure}

\begin{figure}[!tbp]
    \centering
    \includegraphics[width=0.8\linewidth]{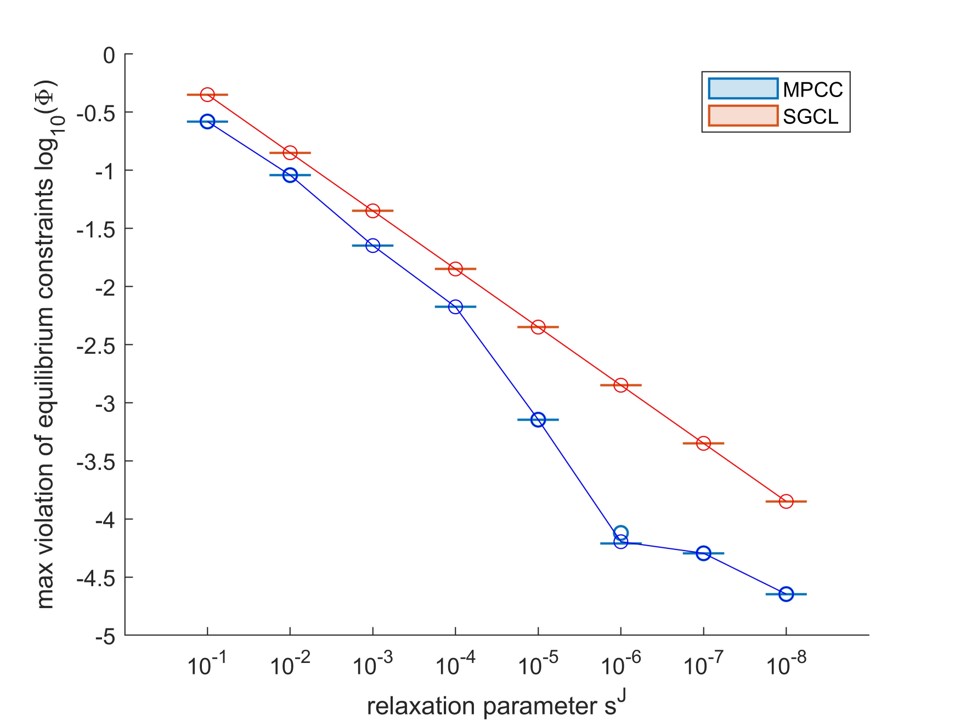}
    \caption{Box plot for relaxation parameter $s^J$ vs. the max violation of equilibrium constraints $\log_{10} \Phi$.}
    \label{fig: param vs equilibrium constraint violation}
\end{figure}

\begin{figure}[!tbp]
    \centering
    \includegraphics[width=0.8\linewidth]{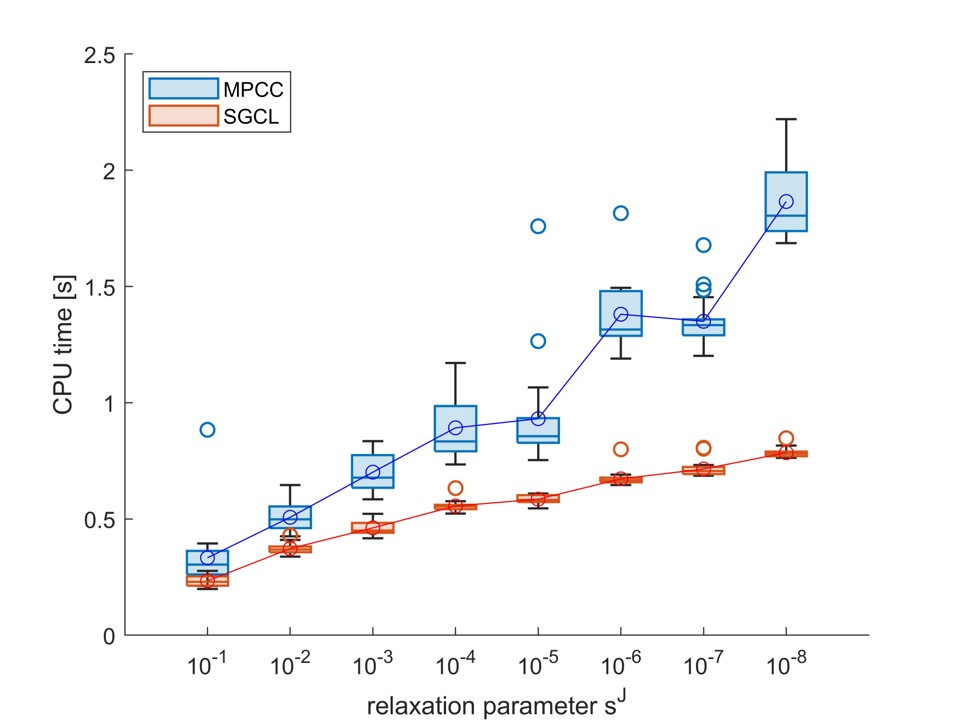}
    \caption{Box plot for relaxation parameter $s^J$ vs. GPU time.}
    \label{fig: param vs time}
\end{figure}

As demonstrated in Figs. \ref{fig: param vs cost} and \ref{fig: param vs equilibrium constraint violation}, at the same $s^J$, both methods find a solution with a nearly identical optimal cost, while the MPCC method has a better equilibrium constraint satisfaction.
Fig. \ref{fig: param vs time} illustrates that the required computation time grows as $s^J$ decreases, and the anomaly of the MPCC method at $s^J = 10^{-6}$ is caused by the frequent calls to the feasibility restoration phase in IPOPT.
Moreover, the MPCC method demands significantly more time than the SGCL method at the same $s^J$.
Since the relaxed problems solved by these two methods exhibit entirely different relaxations for the discretized OCPEC (\ref{equation: discretized OCPEC}), comparing their computation times under the same constraint violations provides a fairer assessment.
In this scenario, the SGCL method (at $s^J = 10^{-7}$) still requires less computation time compared to the MPCC method (at $s^J = 10^{-5}$).
Intuitive interpretations for the effectiveness of the SGCL method are the connectedness of the relaxed feasible set as illustrated in Fig. \ref{fig: regularized gap function region}, and a more concise constraint system as discussed in Subsection \ref{subsection: gap constraint based reformulation for OCPEC}.

\section{Conclusion and future work}\label{section: conclusion and future work}
To solve the OCPEC efficiently, we propose a novel gap-constraint-based reformulation for the OCPEC and solve the discretized OCPEC using an SQP-type method called the SGCL method.
The numerical experiments demonstrate that the SGCL method outperforms the MPCC-tailored solution method in terms of computational efficiency.
The future work will focus on two aspects. 
On the theoretical front, we will analyze the MPCC-tailored stationarity property of the limit point found by the proposed method.
On the practical side, in addition to the discussion in Remark \ref{remark: convex approximation}, we will attempt to design a remedy mode to overcome a common numerical difficulty in SQP-type methods, that is, the infeasible QP subproblem due to the inconsistent constraint linearization.

\bibliography{reference}             
                                                     
\end{document}